\documentclass[a4paper,12pt]{article}

\usepackage{color}

\usepackage{makeidx}

\usepackage{latexsym}
\usepackage{amscd}
\usepackage{amsmath}
\usepackage{amssymb}

\usepackage{amsthm}
\usepackage{float}
\usepackage{graphicx}
\usepackage{psfrag}

\theoremstyle{plain}
\newtheorem{theorem}{Theorem}

\newtheorem{lemma}[theorem]{Lemma}

\newtheorem{proposition}[theorem]{Proposition}

\newtheorem{example}[theorem]{Example}

\theoremstyle{definition}
\newtheorem{definition}[theorem]{Definition}
\newtheorem{remark}[theorem]{Remark}

\newdimen\argwidth
\def\db[#1\db]{%
 \setbox0=\hbox{$#1$}\argwidth=\wd0
 \setbox0=\hbox{$\left[\box0\right]$}
  \advance\argwidth by -\wd0
 \left[\kern.3\argwidth\box0 \kern.3\argwidth\right]}

\newcommand{\bC}{\ensuremath{\mathbb{C}}}

\newcommand{\bN}{\ensuremath{\mathbb{N}}}

\newcommand{\bQ}{\ensuremath{\mathbb{Q}}}
\newcommand{\bR}{\ensuremath{\mathbb{R}}}

\newcommand{\bZ}{\ensuremath{\mathbb{Z}}}

\newcommand{\scE}{\ensuremath{\mathcal{E}}}

\newcommand{\scL}{\ensuremath{\mathcal{L}}}
\newcommand{\scM}{\ensuremath{\mathcal{M}}}

\newcommand{\scP}{\ensuremath{\mathcal{P}}}

\newcommand{\msubt}{\uuu}
\newcommand{\sgn}{\operatorname{sgn}}

\newcommand{\vol}{\operatorname{vol}}
\newcommand{\dist}{\operatorname{dist}}

\newcommand{\const}{{\sf (Const)}}

\newcommand{\elem}{{\sf (Elem)}}
\newcommand{\lin}{{\sf (Lin)}}
\newcommand{\eps}{{\varepsilon}}

\newcommand{\uuu}{\makebox[2ex][c]{\makebox[0pt][c]{$-$}
    \raisebox{.42ex}{\makebox[0pt][c]{$\cdot$}}}}

\title{Periods and elementary real numbers}

\author{Masahiko Yoshinaga\thanks{
Department of Mathematics, 
Graduate School of Science, 
Kobe University, 
1-1, Rokkodai, Nada-ku,
Kobe 657-8501, Japan, 
email: myoshina@math.kobe-u.ac.jp
}}

\date{\today}
\begin{document}

\maketitle

\begin{abstract}
The periods, introduced by Kontsevich and 
Zagier, form a class of complex numbers which contains all 
algebraic numbers and several transcendental 
quantities. Little has been known about 
qualitative properties of periods. In this 
paper, we compare the periods with 
hierarchy of real numbers induced from 
computational complexities. 
In particular we prove that periods 
can be effectively 
approximated by elementary rational Cauchy 
sequences. 
As an application, we exhibit a 
computable real number which is not a period. 
\end{abstract}

\section{Introduction}

In their paper \cite{kz-per}, Kontsevich and Zagier 
introduced the notion of periods: 

\begin{definition}
\label{def:period}
A {\em period} is a complex number whose 
real and imaginary parts are values of absolutely convergent 
integral of rational functions with rational coefficients, 
over domains in $\bR^\ell$ given by polynomial inequalities 
with rational coefficients. 
\end{definition}

The set of all periods is denoted by $\scP\subset\bC$. 
Obviously, $\scP$ is a countable 
set, forms a $\bQ$-algebra (because of Fubini's theorem) 
and contains all algebraic numbers and 
several transcendental quantities, 
like $\pi$ and $\log n$. 
One of their motivations to introduce this notion 
is that the structure of $\scP$ is directly 
related to profound theory of motives. 
See \cite{wald-trans} for related problems in 
transcendental number theory. 

Kontsevich and Zagier pose several conjectures 
and problems on $\scP$. 
However it seems that the qualitative properties 
of $\scP$ have not been well studied so far. 
For instance, they pose the following 
``{\bf Problem 3} {\it Exhibit at least one number which does 
not belong to $\scP$}''. 
We have not had any properties on real numbers 
which can distinguish non-periods from periods.

The purpose of this paper is to give an answer 
to this problem by 
constructing a computable 
real number which can not be a period.

We approach the problem as follows. 
Since the real 
number field $\bR$ is the completion of $\bQ$ with 
respect to the Euclidean norm, 
a positive real number $\alpha\in\bR_{>0}$ can be 
expressed 
as the limit of a positive rational Cauchy sequence 
\begin{equation}
\label{eq:seq}
\lim_{n\rightarrow\infty}\frac{a(n)}{b(n)}=\alpha, 
\end{equation}
where $a$ and $b$ are functions $\bN\rightarrow\bN$. 
Therefore a positive real number $\alpha$ is 
expressed by a pair of functions $a, b:\bN\rightarrow\bN$. 

The observation that 
not all functions $\bN\rightarrow\bN$ 
are computable ``by finite means'', 
since the set of all functions $\bN^\bN$ is 
uncountable, leads us to consider the computability 
of the functions $a$ and $b$. 
The idea of computability goes back to the seminal paper 
\cite{tur-comput} by A. Turing. 
Turing defines computable real numbers as those 
real numbers with computable decimal expansions. 
An equivalent definition is that the real numbers 
which are limits of Cauchy sequences 
(\ref{eq:seq}) with computable functions $a$ and $b$ 
(see \cite{pel-ric} or 
\S\ref{subsec:cr} below). 
So, refined notions of computability enable us to 
hierarchize computable real numbers 
\cite{spe-nicht, rice-recreal, csz-prreal1, csz-prreal2}.

In this paper, we will focus on a proper sub-class called 
``elementary functions'' $\bN\rightarrow\bN$ introduced 
in \cite{csi-elem,kal-elem}  
(see \S\ref{sec:elemreal} below for definitions). 
The main result (Theorem \ref{thm:main}) 
states that every real period is an elementary real 
number, that (roughly speaking) is, 
we can choose $a$ and $b$ from elementary functions. 
And we will also construct a computable 
real number which is not elementary 
(\S\ref{subsec:nonelem}). 
The non-elementary real numbers can not be 
periods by our main result.

Let us briefly describe the idea of the proof. 
First we show that 
periods are generated by the volumes $\vol(D)$ 
of the bounded domains of the form 
$$
D=\left\{(x_1, \ldots, x_\ell)\in\bR^\ell \mid 
G_k(x_1, \ldots, x_\ell)>0, k=1, \ldots, q\right\}, 
$$
where $G_k\in\bZ[x_1, \dots, x_\ell]$ are polynomials 
of integer coefficients. To approximate 
the volume $\vol(D)$, we use the Riemann sum, that is, 
consider the union of small cubes 
$$
V_n:=\mbox{Union of cubes contained in $D$ 
with vertices in }\left(\frac{1}{n}\bZ\right)^\ell. 
$$
Then, clearly, $\vol(V_n)$ converges to $\vol(D)$ as 
$n\rightarrow\infty$. 
However there are two major problems here. 
\begin{itemize}
\item[(a)] Which small cubes are contained 
in the domain $D$? 
\item[(b)] In which rate $\vol(V_n)$ converges to $\vol(D)$? 
(As will be seen in Definition \ref{def:elem1}, 
we have to know the rate of convergence elementarily.) 
\end{itemize}
Let $C\subset\bR^\ell$ be a cube. Then the problem (a) above 
is to ask whether or not the first-order formula 
$$
\forall x(x\in C\Longrightarrow x\in D)
$$
is true. 
In general, the truth assignment for a 
first-order formula with 
quantifiers ($\forall, \exists$) is difficult. 
However, in our situation, Tarski's quantifier elimination for 
real closed ordered field tells us that 
the validity of the above formula can be decided 
by a quantifier free formula. It is simply 
a Boolean combination of 
polynomial inequalities on the coefficients of 
$G_k$'s. 
This enables us to conclude the rational sequence 
$\vol(V_n)$ is elementary. 

The other problem (b) is related to 
count how many small cubes are there near the 
boundary $\partial D$? It is essentially done by 
bound the Minkowski dimension of the 
boundary $\partial D$ by using 
resolution of singularities of algebraic varieties. 

\medskip

The organization of this paper is as follows. 
\S\ref{sec:elemreal} is about elementary functions 
and elementary real numbers. 
Section \S\ref{subsec:cr} begins with the 
definition of the class $\bR_\scE$ 
of real numbers 
computable by a given class $\scE\subset\bN^\bN$ 
of functions. Section \S\ref{subsec:elemf} 
gives the precise definition of elementary 
functions and elementary real numbers. 
In \S\ref{subsec:nonelem}, we algorithmically 
enumerate all elementary Cauchy sequence. 
Then by the diagonal argument, we construct a computable 
real number which is not an elementary 
real number. In view of the main result in 
\S\ref{sec:main}, this number can not be 
a period. In \S\ref{sec:main}, we first state the 
main result. 
After stating the main result 
in \S\ref{subsec:main}, we will reduce the problem 
to the bounded cases by employing results from 
structure theorems of semi-algebraic sets in 
\S\ref{subsec:ur}. In \S\ref{subsec:qe}, we recall 
quantifier elimination by Tarski, and by using it, 
we will construct an elementary rational sequence 
converging to the volume of bounded semi-algebraic domain. 
In the rest, \S\ref{subsec:mc} and \S\ref{subsec:pf}, 
we prove that the sequence converges elementarily.

\section{Elementary real numbers}
\label{sec:elemreal}

{\bf Notation.} In this section, 
$\bN=\{0, 1, 2, 3, \ldots\}$ denotes the 
set of nonnegative integers and 
$(\bN^n)^\bN=\{f:\bN^n\rightarrow\bN\}$ denotes the 
set of all functions from $\bN^n$ to $\bN$. 
We only deal 
with nonnegative real and rational numbers. 

\subsection{Computable real numbers}
\label{subsec:cr}

The set $\bR$ of real numbers is 
defined as the completion of the rational 
number field $\bQ$ by the metric $d(x, y)=|x-y|$. 
In other words, exhibiting a real number is 
equivalent to exhibit a Cauchy sequence 
in $\bQ$. Hence for a given nonnegative real 
number $\alpha\in\bR$, there exist two functions 
$a, b:\bN\rightarrow\bN$ such that 
$$
\lim_{n\rightarrow\infty} 
\frac{a(n)}{b(n)+1}=\alpha. 
$$
(The term ``$+1$'' in the denominator is 
just for avoiding to be equal to zero.) 
In the paper \cite{tur-comput}, Turing introduced 
the notion of computable real numbers by restricting 
the class of functions $a, b:\bN\rightarrow\bN$. 
Following Turing and subsequent studies 
\cite{rose-subrec, csz-prreal1, csz-prreal2}, 
we shall set the following definitions.

\begin{definition}
\label{def:real}
Let $\scE\subset\bN^\bN$ be a class of functions. 
A nonnegative real number $\alpha\in\bR$ is said 
to be $\scE$-computable if 
there exist $a(x), b(x), c(x)\in\scE$ such that 
\begin{equation}
\label{eq:real}
\left|
\frac{a(x)}{b(x)+1}-\alpha
\right|<
\frac{1}{k},\ \mbox{\normalfont for }
\forall x\geq c(k). 
\end{equation}
Denote the set of all 
$\scE$-computable real numbers by $\bR_\scE$. 
\end{definition}

\begin{example}
Obviously $\bR_\scE\subset\bR$ depends on the class $\scE$. \\
{\normalfont (1)} 
Let $\const\subset\bN^\bN$ be the set of 
all constant functions. Then $\bR_\const=\bQ$. 

\medskip

\noindent
{\normalfont (2)} 
Let $\lin\subset\bN^\bN$ be the set of 
all functions of linear growth, that is, 
$$
\lin=\{f\in\bN\mid \exists C>0, \mbox{ s.t. } 
f(n)<C\cdot n\}. 
$$
Then 
$\bR_\lin=\bR$. Indeed for given $\alpha\in\bR$, define 
\begin{eqnarray*}
a(n)&=&\lfloor (n+1)\cdot\alpha\rfloor\\
b(n)&=&n, 
\end{eqnarray*}
which are of linear growth. 
It is easily shown that 
$$
\left|
\frac{a(n)}{b(n)+1}-\alpha
\right|<
\frac{1}{n+1}. 
$$

\medskip

\noindent
{\normalfont (3)} If $\scE$ is the set of all computable 
or recursive 
(resp. primitive recursive) functions, then $\bR_\scE$ is 
the set of computable (resp. primitive recursive) 
real numbers. (See \cite{tur-comput} and \cite{pel-ric}, 
\cite{rice-recreal} for computable numbers. 
And see \cite{csz-prreal1} for a recent survey on 
primitive recursive real numbers.) 
\end{example}

\subsection{Elementary functions}
\label{subsec:elemf}

In order to state the main result, we need the notion 
of {\em elementary functions} 
$\elem$. Here we consider functions 
having any number of arguments, that is, 
$f:\bN^n\rightarrow\bN$ for $n=1, 2, \ldots$. 

We begin with the simplest functions and operations 
on functions. 

\begin{definition}
The zero function: $o(x)=0$. 
The successor function: $s(x)=x+1$. 
The $i$-th projection function: 
$P^n_i(x_1, \ldots, x_n)=x_i$. 
These three functions are called the {\em initial functions}. 
\end{definition}

\begin{definition}
Define the modified subtraction $m:\bN^2\rightarrow\bN$ 
as follows: 
$$
m(x,y)=x\msubt y:=
\left\{
\begin{array}{cc}
x-y& \mbox{if } x\geq y, \\
0  & \mbox{if } x<y. 
\end{array}
\right.
$$
\end{definition}

Let $f(x_1, \ldots, x_m)$ be a function with $m$ arguments. 
Let $g_i(y_1, \ldots, y_n)$ ($i=1, \ldots, m$) be functions 
of $n$ arguments. Then the composition 
$$
f(g_1(y_1, \ldots, y_n), \ldots, 
g_m(y_1, \ldots, y_n))
$$
is a function with $n$ arguments. 

Let $f(t, x_1, x_2, \ldots, x_n)$ be a function with 
$(n+1)$ arguments. We define bounded summation by 
$$
\sum_{t\leq x}f(t, x_1, \ldots, x_n)=
f(0, x_1, \ldots, x_n)+\cdots+ 
f(x, x_1, \ldots, x_n), 
$$
and bounded product by 
$$
\prod_{t\leq x}f(t, x_1, \ldots, x_n)=
f(0, x_1, \ldots, x_n)\times\cdots\times 
f(x, x_1, \ldots, x_n), 
$$
which are functions with $(n+1)$ arguments.

\begin{definition}
The class $\elem$ of elementary functions is the 
smallest class of functions: 
\begin{itemize}
\item[(1)] containing the initial functions, the addition 
$x+y$, the multiplication $x\cdot y$, 
the modified subtraction $x\msubt y$,  
\item[(2)] closed under composition, and 
\item[(3)] closed under bounded summation and product. 
\end{itemize}
\end{definition}

\begin{example}
\label{ex:elem}
The following are examples of elementary functions. 
\begin{itemize}
\item[$(1)$] By definitions, $s(o(x))=1, s(s(o(x)))=2$, etc, 
are elementary. Hence the constant function is elementary. 
Since $s(x)\msubt s(0)=x$, the identity 
function is elementary. The power 
$x^{y+1}=\prod_{k=0}^yx$ is also elementary. 
\item[$(2)$] The sign function 
$$
\sgn(x)=
\left\{
\begin{array}{cc}
1& \mbox{if } x\neq 0, \\
0  & \mbox{if } x=0
\end{array}
\right.
$$
is elementary. Indeed, 
$\sgn(x)=1\msubt(1\msubt x)$. 
\item[$(3)$] Recall that a subset $P\subset\bN^n$ 
is called a predicate. A predicate $P$ is said 
to be elementary if the characteristic function 
$$
\chi_P(x)=
\left\{
\begin{array}{cc}
1& \mbox{if } x\in P, \\
0  & \mbox{if } x\notin P
\end{array}
\right.
$$
is an elementary function. If $P$ and $Q$ are elementary 
predicates, then the Boolean connection $P\wedge Q$, 
$P\vee Q$ and $\neg P$ are also elementary predicates. 
\item[$(4)$] The order predicate 
$$
f_>(x,y)=
\left\{
\begin{array}{cc}
1& \mbox{if } x>y, \\
0  & \mbox{if } x\leq y
\end{array}
\right.
$$
is elementary. Indeed $f_>(x,y)=\sgn(x\msubt y)$. 
Other functions $f_\geq, f_<, f_\leq$ are similarly 
elementary. 
\item[$(5)$] The quotient 
$q(x,y)=\left\lfloor\frac{x}{y+1}\right\rfloor$ 
is elementary. Indeed, 
$$
q(x,y)=\left(\sum_{i=0}^x f_\geq(x, i\cdot (y+1))\right)
\msubt 1. 
$$
Similarly, the logarithm 
$l(a,b)=\lfloor\log_a b\rfloor$ 
and the square root 
$\lfloor\sqrt{x}\rfloor$ are also elementary. 
\item[$(6)$] Bounded minimizer 
$$
(\mu y_1\leq n)(f(y_1, y_2, \ldots, y_k)=0) 
$$
is defined as the least $t\leq n$ such that 
$f(t, y_2, \dots, y_k)=0$ and $n$ if no such $t$. 
If $f:\bN^k\rightarrow\bN$ is elementary, then 
$$
g(n, y_2, \dots, y_k)=(\mu y_1\leq n)(f(y_1, y_2, \ldots, y_k)=0) 
$$
is also elementary. 
\item[$(7)$] The pairing function $J(x,y)$ is defined by 
$$
J(x, y)=\frac{(x+y)(x+y+1)}{2}+y. 
$$
The inverse pairing functions $L(z)$, $R(z)$ 
are defined by the following relations 
$$
J(L(z), R(z))=z,\ L(J(x, y))=x,\ R(J(x,y))=y. 
$$
The functions $L, R$ are also elementary. 
\end{itemize}
\end{example}

\begin{remark}
There also exists a computable but 
non-elementary function, e.g., 
$$
f(x)=x^{x^{\cdot^{\cdot^{\cdot^x}}}} \mbox{ ($x$ floors)}, 
$$
i.e., $f(2)=2^2=4$, $f(3)=3^{3^3}=3^{27}=7625597484987$, 
$f(4)=4^{4^{4^4}}=4^{4^{256}}>4^{1.34078\times 10^{154}}$. 
This is a very rapidly growing function, 
faster than any elementary 
function with one variable. 
(See \cite{rose-subrec} for details.)  
\end{remark}

Recall that the set of elementary real numbers 
$\bR_\elem$ is defined as follows. 

\begin{definition}
\label{def:elem1}
A real number $\alpha\in\bR$ is called elementary if 
there exist $a(x), b(x), c(x)\in\elem$ such that 
\begin{equation}
\label{eq:elem1}
\left|
\frac{a(x)}{b(x)+1}-\alpha
\right|<
\frac{1}{k}, \mbox{\normalfont for }
\forall x\geq c(k). 
\end{equation}
\end{definition}
The following proposition is straightforward. 
\begin{proposition}
The set of elementary real numbers $\bR_\elem$ forms a field. 
\end{proposition}

\begin{definition}
A map $g:\bN\rightarrow\bQ$ is said to be elementary if 
$g$ is expressed as 
$$
g(x)=\frac{a(x)}{b(x)+1}
$$
for some $a(x), b(x)\in\elem$. 
An map $g:\bN\rightarrow\bQ$ is said to be 
{\em fast} if it satisfies 
$$
|g(x)-g(x+1)|<\frac{1}{7^{x+1}},
$$
for $\forall x\in\bN$. 
\end{definition}
\begin{lemma}
\label{lem:fast}
A real number $\alpha\in\bR$ is elementary if and only if 
there exist an elementary fast map 
$g:\bN\rightarrow\bQ$ such that 
\begin{equation}
\label{eq:fast}
\lim_{x\rightarrow\infty}
g(x)=\alpha. 
\end{equation}
\end{lemma}
\proof
Suppose we have $a(x), b(x), c(x)\in\elem$ 
satisfying Eq.(\ref{eq:elem1}). Set $k=8^{n+1}$ and 
$x=c(8^{n+1})$, we have 
$$
\left|
\frac{a(c(8^{n+1}))}{b(c(8^{n+1}))+1}-\alpha
\right|<
\frac{1}{8^{n+1}}. 
$$
Since $a(c(8^{n+1})), b(c(8^{n+1}))$ are elementary on $n$, 
we have $\alpha\in\bR_\elem$. 
Put $g(x)=a(c(8^{x+1}))/(b(c(8^{x+1}))+1)$. Then 
$
|g(n)-\alpha|<8^{-n-1}. 
$
Hence $|g(n)-g(n+1)|<8^{-n-1}+8^{-n-2}$, 
which is less than $7^{-n-1}$. 
\qed

\subsection{A non-elementary real number}
\label{subsec:nonelem}

In this section, we construct a non-elementary 
real number, essentially, by the diagonal argument. 
Together with the main result in the next section, 
it is an example of real number which is not a period. 

First we recall a simpler description of 
elementary functions, due to Mazzanti. 

\begin{proposition}(Mazzanti \cite{mazz-base}) 
All elementary functions can be generated from 
the following four functions by composition: 
\begin{itemize}
\item The successor, $x\mapsto S(x)=x+1$. 
\item The modified subtraction, $(x,y)\mapsto x\msubt y$. 
\item The quotient, 
$(x, y)\mapsto \left\lfloor\frac{x}{y+1}\right\rfloor$. 
\item The exponential function, $(x, y)\mapsto x^y$. 
\end{itemize}
\end{proposition}
Next we enumerate all elementary functions 
$\{f:\bN\rightarrow\bN\mid\mbox{elementary}\}$ of one variable 
by using the pairing functions $J, L, R$ in 
Example \ref{ex:elem} (7). 
For each $e\in\bN$ we attach an elementary function 
$f_e:\bN\rightarrow\bN$ as follows. 
\begin{itemize}
\item[(0)] If $L(e)=0$, then $f_e(x)=x$. (That is, 
$f_{J(0,k)}(x)=x$). 
\item[(1)] If $e=J(1, k)$, then $f_e(x)=S(f_k(x))=f_k(x)+1$. 
\item[(2)] If $e=J(2, k)$, then $f_e(x)=f_{L(k)}\msubt f_{R(k)}$. 
\item[(3)] If $e=J(3, k)$, then 
$f_e(x)=\left\lfloor\frac{f_{L(k)}}{f_{R(k)}+1}\right\rfloor$. 
\item[(4)] If $e=J(4, k)$, then $f_e(x)=(f_{L(k)})^{f_{R(k)}}$. 
\item[(5)] If $e=J(c, k)$ with $c\geq 5$, then $f_e(x)=0$. 
\end{itemize}
\begin{example}
\label{ex:enum}
Here are some examples. 
$f_0(x)=x$ by (0). 
Since $1=J(1,0)$, $f_1=S\circ f_0(x)=x+1$ by (1). 
Since $2=J(0,1)$, again $f_2=x$. 
Since $3=J(2,0)=J(2, J(0,0))$, $f_3=f_0\msubt f_0=0$. 
Since $4=J(1,1)$, $f_4=S\circ f_1=f_1+1=x+2$. 
Since $169=J(1,16)=J(1, J(4,1))=J(1,J(4,J(1,0)))$, 
$f_{169}=f_{16}+1=(x+1)^x+1$, etc. 
\end{example}

Now we can enumerate all elementary maps 
$$
g:\bN\longrightarrow\bQ_{\geq 0}, 
$$
by the following way: 
\begin{equation}
g_e(n):=
\frac{f_{L(e)}(n)}{f_{R(e)}(n)+1}. 
\end{equation}

Obviously the sequence $\{g_e(x)\}_{x\in\bN}$ is not 
Cauchy in general. We enforce being fast on 
these sequences. 
For an elementary sequence 
$g:\bN\rightarrow\bQ$, define 
$\overline{g}$ is 
$$
\overline{g}(n)=
\left\{
\begin{array}{cl}
g(n)&\mbox{ if }(\forall i<n)(|g(i)-g(i+1)|<7^{-i-1}), \\
g(n_0)&\mbox{ otherwise, where } n_0:=
(\mu i<n)(|g(i)-g(i+1)|\geq 7^{-i-1}). 
\end{array}
\right.
$$
The map $\overline{g}:\bN\rightarrow\bQ$ is 
a fast elementary map by definition, and $g$ is 
fast if and only if $g=\overline{g}$. 
\begin{definition}
For $e\in\bN$, define the $e$-th elementary real 
number by 
$$
\beta_e:=\lim_{n\rightarrow\infty}
\overline{g_e}(n). 
$$
\end{definition}
From Lemma \ref{lem:fast}, every elementary real 
number is the limit of a fast sequence, we have 
$$
\{\beta_0, \beta_1, \ldots, \beta_e, \ldots\}=\bR_\elem.  
$$

\begin{example}
\label{ex:40}
First several terms are $\beta_0=0, \beta_1=1, 
\beta_2=\beta_3=0, \beta_4=1/2$ etc. 
Let us compute $\beta_{40}$. Since $40=J(4,4)$, 
$L(40)=R(40)=4$. Thus $g_{40}=f_{40}/(f_{40}+1)$. 
Recall Example \ref{ex:enum} that we already 
have $f_4(x)=x+2$. Hence $g_{40}(x)=\frac{x+2}{x+3}$. 
This is not fast, the enforced one is 
$$
\overline{g_{40}}(x)=
\left\{
\begin{array}{cc}
2/3& \mbox{ if }x=0, \\
3/4& \mbox{ if }x>0. 
\end{array}
\right.
$$
At the end we obtain $\beta_{40}=\frac{3}{4}$. 
\end{example}

Now we construct a non-elementary 
computable real number $\alpha\in\bR$ as the limit 
of sequence 
\begin{equation}
\alpha_n=
\frac{2\varepsilon_1}{3^1}+
\frac{2\varepsilon_2}{3^2}+
\frac{2\varepsilon_3}{3^3}+
\cdots +
\frac{2\varepsilon_n}{3^n}, 
\end{equation}
defined as follows. 
Put $\alpha_0=0$ and define $\varepsilon_n (n\geq 1)$ 
inductively as 
\begin{equation}
\varepsilon_{n+1}=
\left\{
\begin{array}{cc}
0&\mbox{ if }\overline{g_{n}}(n)>\alpha_n+\frac{1}{2\cdot 3^n}\\
1&\mbox{ if }\overline{g_{n}}(n)\leq\alpha_n+\frac{1}{2\cdot 3^n}
\end{array}
\right.
\end{equation}

\begin{proposition}
\label{prop:nonelem}
Set $\alpha=\lim_{n\rightarrow\infty}\alpha_n$, then 
$\alpha\notin\bR_\elem$. 
\end{proposition}
\proof
We shall prove $\alpha\neq\beta_e$ for any $e\in\bN$. 
By the definition of $\alpha_n$, 
$$
\alpha\leq\alpha_n+
2(3^{-n-1}+3^{-n-2}+\cdots )=\alpha_n+3^{-n}. 
$$
So we have 
\begin{equation}
\alpha\in[\alpha_n, \alpha_n+3^{-n}], 
\end{equation}
for all $n\in\bN$. 
Since $|\overline{g_e}(n)-\overline{g_e}(n+1)|<7^{-n-1}$, 
\begin{eqnarray*}
|\overline{g_e}(n)-\beta_e|&<&
7^{-n-1}(1+7^{-1}+7^{-2}+\cdots)\\
&=&\frac{1}{7^n\cdot 6}. 
\end{eqnarray*}
Thus we have 
\begin{equation}
\beta_e\in \left(
\overline{g_e}(n)-\frac{1}{6\cdot 7^n},\ 
\overline{g_e}(n)+\frac{1}{6\cdot 7^n}
\right). 
\end{equation}
If $\overline{g_e}(e)\leq\alpha_{e}+2^{-1}3^{-e}$, then 
$\alpha_{e+1}=\alpha_{e}+2\cdot 3^{-e-1}$. Hence 
$$
\alpha\in\left[
\alpha_{e}+\frac{2}{3^{e+1}},\ \alpha_{e}+\frac{3}{3^{e+1}}
\right]. 
$$
\begin{eqnarray*}
\beta_e&<&\overline{g_e}(e)+\frac{1}{6\cdot 7^e}\\
&\leq&\alpha_{e}+\frac{1}{2\cdot3^{e}}+\frac{1}{6\cdot 7^e}\\
&\leq&\alpha_{e}+\frac{1}{2\cdot3^{e}}+\frac{1}{6\cdot 3^e}\\
&=&\alpha_{e}+\frac{2}{3^{e+1}}\\
&=&\alpha_{e+1}\leq\alpha. 
\end{eqnarray*}
In particular, 
$\alpha\neq\beta_e$. 
If 
$\overline{g_e}(e)>\alpha_{e}+2^{-1}3^{-e}$ we can 
prove $\beta_e>\alpha$ similarly. 
In conclusion we have $\alpha\notin\bR_{\elem}$. 
\qed

The first $80$ terms of the sequence $\varepsilon_n$ are 
the following. 
$$
\begin{array}{c|c|c|c|c|c|c|c|c|c|c|c|c|c|c|c|c|}
\hline
n&1&2&3&4&5&6&7&8&9&10&11&12&13&14&15&16\\
\hline
\varepsilon_n&1&0&1&1&1&1&1&1&0&1&0&1&1&0&1&1\\
\hline
\hline
n&17&18&19&20&21&22&23&24&25&26&27&28&29&30&31&32\\
\hline
\varepsilon_n&0&1&1&1&1&1&1&0&1&1&0&1&0&1&1&0\\
\hline
\hline
n&33&34&35&36&37&38&39&40&41&42&43&44&45&46&47&48\\
\hline
\varepsilon_n&1&1&0&1&0&1&1&1&1&1&1&1&1&1&1&0\\
\hline
\hline
n&49&50&51&52&53&54&55&56&57&58&59&60&61&62&63&64\\
\hline
\varepsilon_n&1&1&0&1&1&1&1&1&1&0&0&1&1&0&1&1\\
\hline
\hline
n&65&66&67&68&69&70&71&72&73&74&75&76&77&78&79&80\\
\hline
\varepsilon_n&0&1&0&1&1&0&1&1&1&0&1&1&1&1&1&1\\
\hline
\end{array}
$$
The real number 
$\alpha/2=\sum_{i=1}^\infty3^{-i}\cdot\varepsilon_i$ 
is not elementary. The first $30$ digits are the following. 
\begin{equation}
\label{eq:nonelem}
\frac{\alpha}{2}=0.388832221773824641256243009581\dots
\end{equation}

\section{Periods are elementary}
\label{sec:main}

\subsection{Main result}
\label{subsec:main}

Now we can state the main result. 

\begin{theorem}
\label{thm:main}
Real periods are elementary real numbers, i.e., 
$$
\scP\subset\bR_\elem. 
$$
\end{theorem}
So the real number $\alpha$ constructed 
above (\ref{eq:nonelem}) 
is not a period.

To prove this theorem, we need to show that 
a given absolutely convergent integration is 
an elementary real number. 
First we will reduce the problem to 
the cases of volumes of bounded semi-algebraic domains. 
Namely, in \S\ref{subsec:ur}, we will prove that $\scP$ is 
generated by volumes $\vol(D)$ of bounded semi-algebraic 
open domains $D\subset\bR^\ell$. The proof is based 
on Hironaka's rectilinearization theorem on 
semi-algebraic sets. Another result, uniformization 
theorem of semi-algebraic sets, is also mentioned 
for later purposes. 

Next step is to construct an elementary 
sequence $\vol(V_n)$ converging to the volume 
$\vol(D)$ of a semi-algebraic domain $D$. 
In \S\ref{subsec:qe} and \S\ref{subsec:rs}, 
this is done by using Riemann sum, that is, 
approximating the domain by small cubes. 
The fact that the sequence is elementary 
is proved by using Tarski's quantifier elimination 
theorem. 

Finally, in \S\ref{subsec:mc} and \S\ref{subsec:pf}, 
we will prove the convergence 
$\vol(V_n)\rightarrow\vol(D)$ is elementary. 
The main task is to count small cubes 
within a $\varepsilon$-neighborhood of the 
boundary $\partial D$. It is closely related 
to estimate the Minkowski dimension and 
the Minkowski content of $\partial D$. 
It is done with the help of uniformization 
theorem for semi-algebraic sets. 
This completes the proof that 
$\vol(D)$ is an elementary real number.

\subsection{Uniformization and rectilinearization}
\label{subsec:ur}

In this section, we recall 
uniformization and rectilinearization theorem 
on subanalytic sets by Hironaka. Our 
main references are \cite{hironaka-sub} 
and \cite{bm-semi}. First let us recall the notions 
of semi-algebraic set and basic open semi-algebraic set. 
(See \cite{bcr-real} for details.)

\begin{definition}
A {\em semi-algebraic subset} of $\bR^\ell$ is a 
finite union of subsets of the form 
\begin{equation}
\label{eq:semialg}
\{
x\in\bR^\ell\mid 
F_1(x)= \dots = F_p(x)=0,\  
G_1(x)>0, \dots, G_q(x)>0
\}, 
\end{equation}
where 
$F_j, G_k\in\bR[x_1, \dots, x_\ell]$. 
\end{definition}
A map from 
a semi-algebraic subset $X\subset\bR^p$ to 
a semi-algebraic subset $Y\subset\bR^q$ is 
called semi-algebraic 
if its graph is a semi-algebraic subset of 
$\bR^{p+q}$. 

\begin{definition}
A {\em basic open semi-algebraic subset} of $\bR^\ell$ is a 
set of the form 
\begin{equation}
\label{eq:basic}
\{
x\in\bR^\ell\mid 
G_1(x)>0, \dots, G_q(x)>0
\}, 
\end{equation}
where 
$G_k\in\bR[x_1, \dots, x_\ell]$. 
\end{definition}

\begin{proposition}
{\normalfont 
\cite[Thm. 5.1.]{bm-semi}
}
\label{prop:unif}
Let $X$ be a closed analytic subset of a real analytic 
manifold $M$. 
Then there is a real analytic manifold $N$ (of the 
same dimension as $X$) and a proper real analytic 
map $\varphi:N\rightarrow M$ such that $\varphi(N)=X.$ 
\end{proposition}

\begin{proposition}
{\normalfont 
\cite[(2.4)]{hironaka-sub}
}
\label{prop:recti}
Let $X$ be a real-analytic space countable at 
infinity. Let $A$ be a globally defined 
semi-analytic set in $X$, i.e., there exists a 
finite system of real analytic functions 
$g_{ij}$ and $f_{ij}$ on $X$ such that 
$$
A=\bigcup_i\{
x\in X\mid g_{ij}(x)=0, f_{ij}(x)>0, \forall j
\}. 
$$
Then there exists a real-analytic map 
$\pi:\widehat{X}\rightarrow X$ such that 
\begin{itemize}
\item[(1)] 
$\widehat{X}$ is smooth and $\pi$ is proper surjective,
\item[(2)] 
for every point $\xi\in\widehat{X}$, there exists a 
local coordinate system $(z_1, \dots, z_n)$ of 
$\widehat{X}$ centered at $\xi$ for which we have: 
within some neighborhood of $\xi$ in $\widehat{X}$, 
$\pi^{-1}(X)$ is a union of quadrants with respect to 
$(z_1, \dots, z_n)$, where a quadrant means a 
set defined by a system of relations $z_1\sigma_1 0, 
z_2\sigma_2 0, \dots, z_n\sigma_n 0$ with 
$\sigma_i$ is either ``$=$'', ``$>$'' or ``$<$''. 
\end{itemize}
\end{proposition}
We note that the map $\pi$ above can be taken to be 
a composition of a finite sequence of blowing-ups 
with smooth centers.

The following apparently 
more general description of $\scP$ is equivalent 
to Definition \ref{def:period}. 
\cite[Thm 2.5, Prop 4.2]{bb-per}: 

\begin{proposition}
\label{prop:period2}
The ring $\scP$ is exactly the ring generated by the 
numbers of the form $\int_\Delta\omega$, where 
$X$ is a smooth algebraic variety of dimension $\ell$ 
defined over $\bQ$, $E\subset X$ is a divisor with 
normal crossings, $\omega\in\Omega^\ell(X)$ is a 
top degree algebraic differential form on $X$, and 
$\Delta\subset X$ is a $\ell$-dimensional 
compact real 
semi-algebraic set with $\partial\Delta\subset E$. 
\end{proposition}
In view of Proposition \ref{prop:recti}, 
we may assume that the semi-algebraic 
cycle $\Delta$ in Proposition \ref{prop:period2} 
is smooth and locally (analytically) a union of quadrants. 

Now we come to prove that real periods are 
elementary. We first reduce the problem to 
the volumes of bounded semi-algebraic sets. 

\begin{lemma}
\label{lem:bdd}
Periods $\scP$ is generated by 
\begin{equation}
\label{eq:basis}
\left\{
\vol(D)\mid
D\subset\bR^k\ \mbox{\normalfont is bounded basic 
open semi-algebraic set}
\right\}. 
\end{equation}
\end{lemma}
\proof 
We will prove: 
\begin{itemize}
\item[(i)] $\scP$ is generated by 
$$
\left\{
\vol(D)\mid
D\subset\bR^k\ \mbox{\normalfont is bounded 
open semi-algebraic set}
\right\}, 
$$
and 
\item[(ii)] The volumes of open semi-algebraic subsets 
of $\bR^\ell$ are generated by those of basic ones. 
\end{itemize}
The second one (ii) is easy. Indeed, a 
semi-algebraic subset of the form (\ref{eq:semialg}) 
with $p>0$ has measure zero. 
As far as we are interested in volumes, 
we can ignore the measure zero sets. 
We may 
consider an open semi-algebraic subset 
is a disjoint union of basic ones modulo 
measure zero sets. 

Now we prove (i). 
We use the description in 
Proposition \ref{prop:period2} 
with $\Delta$ smooth and locally (analytically) 
isomorphic to a union of quadrants. Fix a semi-algebraic 
triangulation $\Delta=\bigcup_\alpha\Delta_\alpha$ and 
also fix base points $p_\alpha\in\Delta_\alpha$ in each simplex. 
By taking the triangulation small enough, 
we may assume that the orthogonal projections 
\begin{equation}
\pi_\alpha:\Delta_\alpha\longrightarrow T_{p_\alpha}\Delta_\alpha
\end{equation}
induce the isomorphism 
$\pi_\alpha:\Delta_\alpha\stackrel{\cong}{\longrightarrow}
\pi_\alpha(\Delta_\alpha)\subset T_{p_\alpha}\Delta_\alpha$. Then 
the image $K_\alpha:=\pi_\alpha(\Delta_\alpha)$ 
is also a semi-algebraic 
set. Denote the inverse of the projection by 
$\psi_\alpha:K_\alpha\stackrel{\cong}{\longrightarrow}\Delta_\alpha$, 
which is a semi-algebraic $C^\infty$-map. 
Fix a coordinate $(z_1, \dots, z_\ell)$ of 
the affine space $T_{p_\alpha}\Delta_\alpha$. Then the pull-back 
of $\omega$ by $\psi_\alpha$ is of the form 
\begin{equation}
(\psi_\alpha)^*\omega=H_\alpha(z)dz_1\wedge\dots\wedge dz_\ell. 
\end{equation}
Since composition and differentiations of 
semi-algebraic functions are also semi-algebraic, 
$H(z)$ is a semi-algebraic $C^\infty$-function. 
So the integration 
$\int_{\Delta_i}\omega=\int_{K_\alpha}H(z)dz_1\dots dz_\ell$ 
is equal to the volume $\vol(D_\alpha)$ of the bounded 
semi-algebraic domain 
\begin{equation}
D_\alpha=\{
(x, t)\in\bR^\ell\times\bR\mid
x\in K_\alpha,\ 0\leq t\leq H_\alpha(x)\}. 
\end{equation}
Thus we have (i). \qed

\subsection{Quantifier elimination}
\label{subsec:qe}

Let $\scL_{OR}$ be the language 
$$
\scL_{OR}=(+, -, \cdot, 0, 1, <, =)
$$
of ordered rings. We consider the theory $T$ of 
real number field $\bR$ with the language $\scL_{OR}$. 
Recall that a quantifier free formula 
$\psi(x_1, \dots, x_n)$ is a Boolean combination of 
inequalities $p(x_1, \dots, x_n)>0$, where 
$p\in\bZ[x_1, \dots, x_n]$. 
The following is due to Tarski \cite{tar-dec}, 
see also \cite{coh-dec}. 

\begin{theorem}
\label{thm:qe}
{\normalfont (Tarski)} 
On the real number field, 
every $\scL_{OR}$-formula $\varphi(x_1, \ldots, x_n)$ is 
equivalent to a quantifier free formula 
$\varphi^*(x_1, \ldots, x_n)$, i.e., 
$$
T\vDash
\forall x_1\forall x_2\cdots\forall x_n
\left(
\varphi(x_1, \ldots, x_n)\Leftrightarrow\varphi^*(x_1, \ldots, x_n)
\right). 
$$\qed
\end{theorem}

Let 
\begin{equation}
D=\{{x}=(x_1, \ldots, x_\ell)\in\bR^\ell\mid 
G_k({x})>0,\ k=1, \ldots, q\}, 
\end{equation}
be a domain in $\bR^\ell$, 
where $G_k(x)\in\bZ[x_1, \ldots, x_\ell]$ and 
set 
$$
G_k(x)=
\sum_{J}a_{kJ}x^J, 
$$
where $J=(j_1, \ldots, j_\ell)$ is multi-index and 
denoting $x^J=x_1^{j_1}\cdots x_\ell^{j_\ell}$. 

Let us consider the next predicates with 
variables $s_i, t_i, a_{kJ}$, 

\medskip

\noindent
$R(s_i, t_i, a_{kJ}: 1\leq i\leq\ell, 1\leq k\leq q, J)$: 
\begin{equation}
\forall x_1\dots\forall x_\ell
\left(
(s_i\leq x_i\leq t_i), 
i=1, \ldots, \ell
\Rightarrow
(G_k(x)>0), 
k=1, \ldots, q
\right)
\end{equation}
The above formula means that the box 
$\prod_{i=1}^\ell[s_i, t_i]$ is 
contained in the domain $D$, 
\begin{equation}
\prod_{i=1}^\ell[s_i, t_i]=
[s_1, t_1]\times\cdots\times [s_\ell, t_\ell]
\subset D. 
\end{equation}

From Theorem \ref{thm:qe}, we have a 
quantifier free formula 
$R^*(s_i, t_i, a_{kJ})$ which satisfies 
for $\forall s_i, t_i, a_{k,J}$, 
\begin{equation}
R^*(s_i, t_i, a_{kJ})\Longleftrightarrow 
[s_1, t_1]\times\cdots\times [s_\ell, t_\ell]
\subset D. 
\end{equation}

\subsection{Riemann sum}
\label{subsec:rs}

Let $D\subset\bR^\ell$ be a basic open 
semi-algebraic subset as in (\ref{eq:basic}). 
Now we assume that $D$ is bounded and 
contained in a large cube $[0,r]^\ell$, $r>0$.

Then the volume $\vol(D)$ is approximated by 
the inner Riemann sum. 
For positive an integer $n>0$ and $k_1, \dots, k_\ell\in\bN$, 
define a small cube $C_n(k_1,\dots, k_\ell)$ of size  
$r/n$ by 
$$
C_n(k_1, \ldots, k_\ell)=
\left[\frac{k_1r}{n}, \frac{(k_1+1)r}{n}\right]\times
\cdots\times
\left[\frac{k_\ell r}{n}, \frac{(k_\ell +1)r}{n}\right]. 
$$
Trivially these cubes subdivides the large cube 
$[0,r]^\ell=\bigcup_{0\leq k_i<n}C_n(k_1,\dots, k_\ell)$. 
Let us denotes by $V_n$ the union 
$$
V_n=\bigcup_{C_n(k)\subset D}C_n(k_1, \dots, k_\ell) 
$$
of small cubes which are contained in $D$. 
We will prove that $\vol(V_n)\rightarrow\vol(D)$ 
(as $n\rightarrow\infty$) determines an elementary 
real number. 
\begin{lemma}
The function 
$$
\bN\longrightarrow\bQ\ \left(
n\longmapsto\vol(V_n)\right)
$$
is elementary. 
\end{lemma}
\proof
To compute Riemann sum $\vol(V_n)$, 
we have to know for which 
$(k_1, \dots, k_\ell)$ the small cube $C_n(k_1, \dots, k_\ell)$ 
is contained in $D$. From Theorem \ref{thm:qe} in the 
previous section, 
this is decided by a quantifier free 
formula $R^*(s_i, t_i, a_{kJ})$. 
By definition, it is a Boolean 
combination of the predicates of the form 
\begin{equation}
p(s_i, t_i, a_{kJ})>0, 
\end{equation}
with $p\in\bZ[s_i, t_i, a_{kJ}]$. 
The truth value of the statement $C_n(k)\subset D$ is decided 
by checking the truth values of 
Boolean combination of predicates of the form 
\begin{equation}
p\left(
\frac{k_ir}{n}, \frac{(k_i+1)r}{n}, a_{kJ}
\right)>0. 
\end{equation}
Thus the relation $C_n(k_1, \dots, k_\ell)\subset D$ 
can be decided elementarily, that is, there exists 
an elementary function 
$$
\varphi:\bN^{\ell+1}\longrightarrow\bN,\ 
\left(
(n, k_1, \dots, k_\ell)\longmapsto 
\varphi(n, k_1, \dots, k_\ell)
\right)
$$
such that 
\begin{equation}
\varphi(n, k_1, \dots, k_\ell)=
\left\{
\begin{array}{cc}
1&\mbox{ if }C_n(k_1, \dots, k_\ell)\subset D, \\
0&\mbox{otherwise}. 
\end{array}
\right.
\end{equation}

Thus the volume $\vol(V_n)$ of the union of small cubes 
is expressed as 
\begin{equation}
\vol(V_n)=\left(\frac{r}{n}\right)^\ell
\sum_{0\leq k_i\leq n}
\varphi(n, k_1, \ldots, k_\ell), 
\end{equation}
which is an elementary function on $n$. 
\qed

Next we have to estimate the rate of convergence 
$$
\lim_{n\rightarrow\infty}\vol(V_n)=\vol(D). 
$$

\subsection{Minkowski content}
\label{subsec:mc}

In this subsection we recall notations on 
Minkowski dimension and Minkowski contents from 
\cite{kp-dom}. 

First let $B=\{x\in\bR^N\mid |x|<1\}$. 
Then 
$$
\Upsilon_N:=\vol(B)=
2\frac{\pi^{N/2}}{N\cdot\Gamma(N/2)}. 
$$
Here $\scL^N$ denotes the $N$-dimensional Lebesgue measure. 

\begin{definition}
Suppose $A\subset\bR^N$ and $0\leq K\leq N$. 
The {\em $K$-dimensional upper Minkowski content of $A$}, 
denoted by $\scM^{*K}(A)$, is defined by 
$$
\scM^{*K}(A)=
\limsup_{\eps\downarrow 0}
\frac{\scL^N\{x\mid\dist(x,A)<\eps\}}{\Upsilon_{N-K}\eps^{N-K}}
$$
\end{definition}

\begin{proposition}
\label{prop:mink}
{\normalfont (\cite[Prop 3.5.5]{kp-dom})} Let 
$f:\bR^K\rightarrow\bR^N$ be a $C^1$-map. 
$A\subset\bR^K$ is compact with 
$$
A\subset\{x\mid |D(f)|\leq\rho\}, 
$$
then 
$$
\scM^{*K}(f(A))\leq\rho^K\scL^K(A). 
$$
\end{proposition}

\subsection{Proof, completion}
\label{subsec:pf}

Now we return to the proof of 
Theorem \ref{thm:main}, $\scP\subset\bR_\elem$. 
In view of Lemma \ref{lem:bdd}, 
it is enough to show that 
the sequence $\vol(V_n)\rightarrow\vol(D)$ constructed 
in \S\ref{subsec:rs} converges effectively. 
The following lemma concludes 
$\vol(D)\in\bR_\elem$. 

\begin{lemma}
\label{lem:conv}
There exists a constant $L=L(D)$ depending only on $D$, 
such that if $k$ and $n$ satisfy 
\begin{equation}
\label{eq:condition}
4rL\sqrt{\ell}k<n, 
\end{equation}
then $|\vol(D)-\vol(V_n)|<1/k$. 
\end{lemma}
\proof 
Set $P(x)=\prod_{k=1}^q G_k(x)$. 
Then $\partial D\subset \{P=0\}$. 
By the Uniformization theorem 
(Proposition \ref{prop:unif}), $X=\{P=0\}$ is an 
image $\pi(X)$ of a proper analytic map 
$\pi:X\rightarrow \bR^\ell$. Since $\partial D\subset \bR^\ell$ 
is compact, from Proposition \ref{prop:mink}, 
the $(\ell-1)$-dimensional Minkowski content 
$\scM^{*(\ell-1)}(\partial D)$ of 
the boundary $\partial D$ is finite. 
There is a constant $L>0$ and $\eps_0>0$ 
such that 
$$
\frac{\scL^{\ell}(\{y\in\bR^\ell\mid 
\dist (y, \partial D)<\eps\})}{2\eps}<L, 
$$
for $0<\eps<\eps_0$. Equivalently we have 
$$
\scL^{\ell}(\{y\in\bR^\ell\mid 
\dist (y, \partial D)<\eps\})<2\eps L. 
$$
Choose $n$ large enough and $\eps$ as 
\begin{equation}
\label{eq:large}
\frac{r\sqrt{\ell}}{n}=\frac{\varepsilon}{2}, 
\end{equation}
note that the LHS is exactly the diagonal length of the 
small cube $C_n(k_1, \dots, k_\ell)$. 
Let us consider the subset of $D$ which is 
$\varepsilon$-away from the 
boundary (or removing $\varepsilon$-neighborhood of 
the boundary) 
\begin{equation}
D_{>\varepsilon}=
\{x\in D\mid 
\dist(x, \partial D)>\varepsilon\}. 
\end{equation}
It is easily seen that, under (\ref{eq:large}), 
\begin{equation}
D_{>\varepsilon}\subset V_n\subset D. 
\end{equation}
Instead of $\vol(D-V_n)$, we will estimate 
$\vol(D-D_{>\eps})$.

Hence if we choose $n$ as in Eq. (\ref{eq:large}), 
\begin{eqnarray*}
|\vol(D)-\vol(V_n)|&<&|\vol(D)-\vol(D_{>\eps})|\\
&=&\scL^{\ell}(\{y\in D \mid \dist (y, \partial D)<\eps\})\\
&<&\scL^{\ell}(\{y\in \bR^\ell \mid \dist (y, \partial D)<\eps\})\\
&<&2\eps L\\
&=&\frac{4r\sqrt{\ell}L}{n}. 
\end{eqnarray*}
Thus if (\ref{eq:condition}) is satisfied, 
we have $|\vol(D)-\vol(V_n)|<1/k$. 
\qed

\medskip

\noindent
{\bf Acknowledgment.} 
The author thanks to 
Professor Masa-Hiko Saito for his interests to 
this work and constant 
encouragements. The author also thanks to 
Professor 
Toshiyasu Arai, 
Makoto Kikuchi, 
Takefumi Kondo, 
Hiraku Kawanoue, 
Takeshi Nozawa, 
Okihiro Sawada 
for comments and useful conversations on 
several topics treated in this paper.

\bibliographystyle{plain}

\begin{thebibliography}{10}

\bibitem{bb-per}
Prakash Belkale and Patrick Brosnan.
\newblock Periods and {I}gusa local zeta functions.
\newblock {\em Int. Math. Res. Not.}, (49):2655--2670, 2003.

\bibitem{bm-semi}
Edward Bierstone and Pierre~D. Milman.
\newblock Semianalytic and subanalytic sets.
\newblock {\em Inst. Hautes \'Etudes Sci. Publ. Math.}, (67):5--42, 1988.

\bibitem{bcr-real}
Jacek Bochnak, Michel Coste, and Marie-Fran{\c{c}}oise Roy.
\newblock {\em Real algebraic geometry}, volume~36 of {\em Ergebnisse der
  Mathematik und ihrer Grenzgebiete (3) 
}.
\newblock Springer-Verlag, Berlin, 1998.
\newblock Translated from the 1987 French original, Revised by the authors.

\bibitem{csz-prreal1}
Qingliang Chen, Kaile Su, and Xizhong Zheng.
\newblock Primitive recursive real numbers.
\newblock {\em MLQ Math. Log. Q.}, 53(4-5):365--380, 2007.

\bibitem{csz-prreal2}
Qingliang Chen, Kaile Su, and Xizhong Zheng.
\newblock Primitive recursiveness of real numbers under different
  representations.
\newblock In {\em Proceedings of the Third International Conference on
  Computability and Complexity in Analysis (CCA 2006)}, volume 167 of {\em
  Electron. Notes Theor. Comput. Sci.}, pages 303--324 (electronic), Amsterdam,
  2007. Elsevier.

\bibitem{coh-dec}
Paul~J. Cohen.
\newblock Decision procedures for real and {$p$}-adic fields.
\newblock {\em Comm. Pure Appl. Math.}, 22:131--151, 1969.

\bibitem{csi-elem}
Paul Csillag.
\newblock Eine {B}emerkung zur {A}ufl\"osung der eingeschachtelten {R}ekursion.
\newblock {\em Acta Univ. Szeged. Sect. Sci. Math.}, 11:169--173, 1947.

\bibitem{hironaka-sub}
Heisuke Hironaka.
\newblock Subanalytic sets.
\newblock In {\em Number theory, algebraic geometry and commutative algebra, in
  honor of Yasuo Akizuki}, pages 453--493. Kinokuniya, Tokyo, 1973.

\bibitem{kal-elem}
L\'aszl\'o Kalm\'ar.
\newblock Egyszer\"u p\'elda eld\"onthetetlen aritmetikai probl\'em\'ara.
\newblock {\em Matematikai \'es Fizikai Lapok}, 50:1--23, 1943.

\bibitem{kz-per}
Maxim Kontsevich and Don Zagier.
\newblock Periods.
\newblock In {\em Mathematics unlimited---2001 and beyond}, pages 771--808.
  Springer, Berlin, 2001.

\bibitem{kp-dom}
Steven~G. Krantz and Harold~R. Parks.
\newblock {\em The geometry of domains in space}.
\newblock Birkh\"auser Advanced Texts: Basler Lehrb\"ucher. 
Birkh\"auser Boston Inc., Boston, MA, 1999.

\bibitem{mazz-base}
Stefano Mazzanti.
\newblock Plain bases for classes of primitive recursive functions.
\newblock {\em MLQ Math. Log. Q.}, 48(1):93--104, 2002.

\bibitem{pel-ric}
Marian~B. Pour-El and J.~Ian Richards.
\newblock {\em Computability in analysis and physics}.
\newblock Perspectives in Mathematical Logic. Springer-Verlag, Berlin, 1989.

\bibitem{rice-recreal}
H.~G. Rice.
\newblock Recursive real numbers.
\newblock {\em Proc. Amer. Math. Soc.}, 5:784--791, 1954.

\bibitem{rose-subrec}
H.~E. Rose.
\newblock {\em Subrecursion: functions and hierarchies}, volume~9 of {\em
  Oxford Logic Guides}.
\newblock The Clarendon Press Oxford University Press, New York, 1984.

\bibitem{spe-nicht}
Ernst Specker.
\newblock Nicht konstruktiv beweisbare {S}\"atze der {A}nalysis.
\newblock {\em J. Symbolic Logic}, 14:145--158, 1949.

\bibitem{tar-dec}
Alfred Tarski.
\newblock {\em A decision method for elementary algebra and geometry}.
\newblock University of California Press, Berkeley and Los Angeles, Calif.,
  1951.
\newblock 2nd ed.

\bibitem{tur-comput}
A.~M. Turing.
\newblock On computable numbers, with an application to the
  "entscheidungsproblem".
\newblock {\em Proc. London Math. Soc.}, 42:230--265, 1936.

\bibitem{wald-trans}
Michel Waldschmidt.
\newblock Transcendence of periods: the state of the art.
\newblock {\em Pure Appl. Math. Q.}, 2(2):435--463, 2006.

\end{thebibliography}


\end{document}